# On multidimensional item response theory – a coordinate free approach

## Tamás Antal

*4 Bayberry Road, Princeton, NJ 08540*
*e-mail:* antaltamas@gmail.com

**Abstract:** A coordinate system free definition of complex structure multidimensional item response theory (MIRT) for dichotomously scored items is presented. The point of view taken emphasizes the possibilities and subtleties of understanding MIRT as a multidimensional extension of the "classical" unidimensional item response theory models. The main theorem of the paper is that every monotonic MIRT model looks the same; they are all trivial extensions of univariate item response theory.

**AMS 2000 subject classifications:** Primary 62H05; secondary 62P25.

**Keywords and phrases:** geometric methods, statistical modeling.

## 1. Introduction

Complex structure multidimensional item response theory (MIRT) is built on the idea that a single item, however simple it might be, carries the possibility of an inner structure. That is, in usual terminology one speculates that it is possible to measure several cognitive areas with one item. The number of cognitive areas so measured may vary among items, even though usual models assume that it is fixed for a collection of items (a test) and let a factor analysis type procedure decide on the number and mixture of cognitive areas measurable by the items.

The point of view taken in this note is that any unidimensional item response theory (IRT) model can be thought of as a specialization of a MIRT model. Hence, the major task is to identify how much of the well established tools and nomenclature of unidimensional IRT can be preserved in the multidimensional context and, from the other direction, how different multidimensional notions may specialize to the same unidimensional entity. When the latter happens, that is when two different multidimensional objects yield the same unidimensional specialization, then both multidimensional notions could be considered proper generalizations of the underlying unidimensional quantity. A careful study should then be devised to decide which generalization is more appropriate with respect to the application at hand.

There is, on the other hand, the possibility of not finding proper multidimensional generalization for some unidimensional notions. This topic also deserves careful research and understanding.

Here, we consider what is termed *complex structure* MIRT. Usually, IRT models have two components: the item likelihood and the population distribution.

In *simple structure* MIRT one item represents only one dimension and without a multivariate population distribution the entire likelihood of the model would factor as a product of univariate pieces. In complex structure MIRT this factorization is impossible, by definition, irrespective of population model chosen. Our main theorem will hold irrespective of the complexity of the structure.

The structure of the paper is as follows. A short overview of unidimensional IRT is followed by the *absolute*, that is coordinate system free, definition of MIRT. The connection with the usual approach is also shown via a discussion of two widely accepted models. Then, the development of the main thesis follows. In this we prove that MIRT models are all alike and they all can be obtained as a trivial extension of an appropriate unidimensional item response theory model. Two sections on some thoughts about capturing cognitive dimensions and on understanding the role of the notion of dimension-wise independence close the presentation.

## 2. Unidimensional Item Response Theory

To make the generalization to the multidimensional framework easier, let us first summarize some features of unidimensional IRT. Measurement takes place during the formation of the response matrix $X \in M_{N \times I}(\mathbb{N})$ with elements $x_{ni} \in \mathbb{N}$ for student $n = 1, \ldots, N$ and item $i = 1, \ldots, I$. In a dichotomous setting (which is assumed throughout the paper to simplify the presentation) $x_{ni} = 1$ if student $n$ responded correctly to item $i$, otherwise it is zero. As a major simplification of the modeling of the cognitive process it is assumed that the response to an item is stochastically determined by the ability $\theta$ and item parameters $\beta_i := (a_i, b_i, c_i)$ via the item response function ([1]):

$$P^{3\mathrm{pl}}(\theta, \beta_i) := \mathrm{Prob}(x_{ni} = 1 \mid \theta, \beta_i) = c_i + \frac{1 - c_i}{1 + e^{-a_i(\theta - b_i)}}. \tag{1}$$

There are, of course, many different item response functions in use, the three parameter logistic model is chosen here only as an illustration. The other substantial simplification used in building the model is the assumption of independence of conditional probabilities $P^{3\mathrm{pl}}_{ni}$ across an arbitrary subset $S \subset \{1, \ldots, N\} \times \{1, \ldots, I\}$ of student-item pairs.

The two most popular models built out of these blocks are the joint unidimensional IRT and the marginal unidimensional IRT. Joint IRT states that the total likelihood depends explicitly on the ability of the given students:

$$L^{\mathrm{joint}}(X; \Theta, \mathcal{B}) = \prod_{n,i} P^{3\mathrm{pl}}(\theta_n, \beta_i)^{x_{ni}} (1 - P^{3\mathrm{pl}}(\theta_n, \beta_i))^{1 - x_{ni}} \tag{2}$$

with corresponding log-likelihood:

$$\mathcal{L}^{\mathrm{joint}}(X; \Theta, \mathcal{B}) = \sum_{n,i} x_{ni} \log(P^{3\mathrm{pl}}(\theta_n, \beta_i)) + (1 - x_{ni}) \log(1 - P^{3\mathrm{pl}}(\theta_n, \beta_i)). \tag{3}$$

Here, $\Theta = (\theta_1, \ldots, \theta_N)$ and $\mathcal{B} = (\beta_1, \ldots, \beta_I)$ are the collections of all abilities and item parameters, respectively.

In the marginal theory the likelihood depends only on the distributional properties of student's population:

$$L^{\text{marg}}(X; \mathcal{B}, \Phi) = \prod_n \int_{\mathbb{R}} \prod_i P^{3\text{pl}}(\theta, \beta_i)^{x_{ni}}(1 - P^{3\text{pl}}(\theta, \beta_i))^{1-x_{ni}} \mathrm{d}\mu_n(\theta) \qquad (4)$$

with log-likelihood

$$\mathcal{L}^{\text{marg}}(X; \mathcal{B}, \Phi) = \sum_n \log \int_{\mathbb{R}} \prod_i P^{3\text{pl}}(\theta, \beta_i)^{x_{ni}}(1 - P^{3\text{pl}}(\theta, \beta_i))^{1-x_{ni}} \mathrm{d}\mu_n(\theta), \quad (5)$$

where $\mu_n$ is the density measure of student $n$ over $\mathbb{R}$ and $\Phi$ is the collection of distributional parameters for student's ability. In parametric setting usually $\mu_n$ is given as

$$d\mu_n(\theta) = \varphi_n(\theta)\mathrm{d}\theta$$

with some density function $\varphi_n$.

The quantities

$$L_n^{\text{st}}(X; \theta, \mathcal{B}) = \prod_i P^{3\text{pl}}(\theta, \beta_i)^{x_{ni}}(1 - P^{3\text{pl}}(\theta, \beta_i))^{1-x_{ni}} \qquad (6)$$

and

$$L_i^{\text{it}}(X; \Theta, \beta) = \prod_n P^{3\text{pl}}(\theta_n, \beta)^{x_{ni}}(1 - P^{3\text{pl}}(\theta_n, \beta))^{1-x_{ni}} \qquad (7)$$

are the student and item likelihoods, respectively.

A maximization of the joint model can be achieved by iteratively maximizing all the student likelihoods with fixed item parameters to obtain the next approximation of the abilities and all the item likelihoods with fix abilities to obtain the next approximation of item parameters. Starting values can be constructed from careful item analysis.

It is worthwhile to analyze the shape of the student likelihood function. It is a product of the conditional probabilities of the actual responses over all the items administered to the student. As a function of $\theta$ the probability of the correct response is increasing when the actual response is correct and decreasing for incorrect actual response. As a consequence, a student likelihood will be increasing if all the actual responses are correct and decreasing if all the actual responses are incorrect. This in turn pushes the location of the maximum likelihood solution for the given student to plus or minus infinity. For the item likelihood a similar statement holds. When at least two responses are different in each row and in each column of the dichotomous response matrix the existence of the unique maximum place is guaranteed in every step of the iteration. This, however does not necessarily imply that the iterative method will be convergent ([4] gives a necessary and sufficient condition for the convergence of the joint Rasch model). The student likelihood can be well approximated by a normal distribution (especially when the number of items is large enough) and its curvature will be inversely proportional to the asymptotic standard error of the ability estimates.

The marginal model does not suffer from the same restriction so severely, because the population density function, when chosen according to usual practice, will be sufficient to ensure the existence of a finite maximum place at each step of the iteration. In this case, the standard errors associated with the given row or column will be higher when constant response pattern is present.

For this discussion to even make sense, we had to use the trivially available ordering of real numbers (playing the role of ability space in the unidimensional case) to use such notion as "increasing ability". This point will be central to the multidimensional extension, since there will be no natural choice of ordering of multidimensional abilities.

## 3. Multidimensional Item Response Theory

### 3.1. Introduction

In what follows we explore the possibility of defining MIRT in geometric terms without direct reference to coordinates. As before, the full response likelihood for a student is formed by multiplying single item conditional probabilities together (invoking the assumption of local independence). The likelihood of an MIRT model is then constructed by incorporating some sort of population model with these response likelihoods similar to the joint and marginal univariate cases (Equations 2, 4).

The classification of MIRT is achieved at the level of a single item conditional probability in the same way as we would characterize a univariate IRT model as Rasch, 2PL or 3PL model. This does not mean that we restrict our presentation to single item *tests*. Realistic tests are treated using the local independence assumption as discussed before (Equations 2, 4, and 6).

With this now clarified, from what follows, unless noted otherwise, we shall drop any reference to any particular student and item. This will also help us avoiding overflow of indices in the multidimensional context.

### 3.2. Basic Models

Even though widely investigated, MIRT is not yet widespread as an operational model. Hence, identifying the major players among the competing MIRT models is difficult. Here, only two models are discussed, one by [12] and another one by [10].

First, for an item we associate a vector of discriminations $a = (a_1, \ldots, a_D) \in \mathbb{R}^D$ and a vector of difficulties $b = (b_1, \ldots, b_D) \in \mathbb{R}^D$. With these the functional representation of the *dimension-wise independent* MIRT response likelihood of

[12] has the form

$$f^w_{a,b} : \mathbb{R}^D \to [0,1], \quad \theta \mapsto f^w_{a,b}(\theta) = \prod_{d=1}^{D} \frac{1}{1+e^{-a_d(\theta_d-b_d)}}, \tag{8}$$

where $\theta = (\theta_1, \ldots, \theta_D) \in \mathbb{R}^D$. If the conditional probability of passing the $d$th dimension of the item is given by $\frac{1}{1+e^{-a_d(\theta_d-b_d)}}$, then (8) can indeed be understood as the joint probability of passing *all* the independent dimensions of the item. Unless there are separate observed scores for each dimension, language like "correct response on dimension $d$" cannot be used. In lack of this we used the "passing a dimension" term, which may refer to an unobservable event.

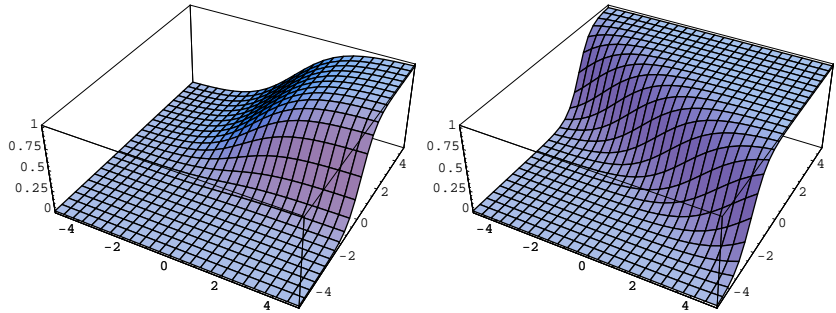

Fig 1. *Dimension-wise independent and Scalar Product MIRT hypersurface*

[8] (see also [10]) put forward a model in which the response likelihood takes the functional representation

$$f^{sp}_{a,b} : \mathbb{R}^D \to [0,1], \quad \theta \mapsto f^{sp}_{a,b}(\theta) = \frac{1}{1+e^{-\langle a \mid \theta \rangle - b}}, \tag{9}$$

where $a$ is as before and $b \in \mathbb{R}$. $\langle x \mid y \rangle = \sum_{d=1}^{D} x_d y_d$ is the usual scalar product of $x, y \in \mathbb{R}^D$. We use the term *Scalar Product MIRT* to refer to this model.

As a last step before embarking on the dimension free definition of MIRT let us write the marginal likelihood of the Scalar Product model assuming multivariate normal population distribution. Using the notation of the previous section, the conditional probability of the response $x_{ni}$ for item $i$ and student $n$ is given by

$$P\left(x_{ni} \mid \theta, \beta_i = (a_i, b_i) \in \mathbb{R}^D \times \mathbb{R}\right) = \frac{1}{1+e^{-(2x_{ni}-1)(\langle a_i \mid \theta \rangle + b_i)}}. \tag{10}$$

Then, the likelihood of the Scalar Product MIRT model is given by

$$L(X; \mathcal{B}, \Phi) = \prod_n \int_{\mathbb{R}^D} \prod_i P(x_{ni} \mid \theta, \beta_i) \varphi(\theta; \nu_n, \Sigma_n) \mathrm{d}^D \theta, \tag{11}$$

where $\varphi(\theta; \nu_n, \Sigma_n)$ is the multivariate normal density with (possibly) student dependent mean $\nu_n$ and covariance $\Sigma_n$. Moreover, $\mathcal{B} = (\beta_1, \ldots, \beta_I)$ is the collection of item parameters and $\Phi = (\nu_1, \Sigma_1, \ldots, \nu_N, \Sigma_N)$ is that of the population parameters.

The multidimensional student likelihood is given by the product of individual conditional probabilities over items administered to a given student $n$:

$$\prod_i P(x_{ni} \mid \theta, \beta_i). \tag{12}$$

Figure 4 depicts a possible student likelihood in the two dimensional case.

The likelihood (11) is a multidimensional generalization of the univariate marginal likelihood given by (4).

### *3.3. Definition of MIRT*

Our goal in this section is to give a definition of MIRT with as few assumption as possible. Multidimensional item response theory postulates that with one single item multiple cognitive abilities could be detected. To accommodate this idea, one has to change the model for the ability space from the one dimensional vector space $\mathbb{R}$ to a finite dimensional vector space $V_\theta$. While any finite dimensional vector space $V$ is linearly isomorphic to $\mathbb{R}^D$ for $D = \dim(V)$ (see (14) for an explicit way of constructing such an isomorphism), this isomorphism is not canonical (there is not a unique isomorphism $V \to \mathbb{R}^D$). By this, and other reasons that will become clear as we proceed, we chose not to use $\mathbb{R}^D$ as a mathematical model of ability space.

The reader unfamiliar with these notions is referred to [5] for an excellent introduction to linear algebra. Also, an intuitive understanding of the basic notions of smooth manifolds should help understanding of what follows, although not strictly necessary. Among the many fine references to the topic the interested reader may find [11] useful.

The basic object in unidimensional IRT is the item response function (IRF) and its graph, the item response curve (IRC). Recall, that the graph of a function $f : A \to B$ is a subset $\text{graph}(f) = \{(x, f(x)) \in A \times B \mid x \in A\}$ of $A \times B$. IRC is a one dimensional smooth submanifold of $\mathbb{R} \times [0, 1]$. While there is a scaling freedom even in the one dimensional case (e.g. the (in)famous 1.7 multiplier in the logistic models), the possibility of ambiguous interpretation is minimal and one may use the functional (IRF) and the geometrical (IRC) representation almost interchangeably.

In the multidimensional case, however, the matter is not so straightforward. As we shall see, the functional and the geometric representations are different in a subtle way. One way to keep the presentation coordinate system free in multidimensional IRT is to postulate that the theory is given by an *item response hypersurface (IRHS)*. As in the unidimensional case, the IRHS is used to express the probability of correct response given an ability in $V_\theta$.

Before defining this notion, let us fix some notations. For any $v \in V_\theta$ the *ray* of $v$ is defined to be the line $\mathbb{R} \cdot v$ in $V_\theta$ determined by $v$: $\mathbb{R} \cdot v = \{\lambda v \in V_\theta \mid \lambda \in \mathbb{R}\}$. Similarly, for $v, w \in V_\theta$ the $v$-directed line going through $w$ is defined by

$$w + \mathbb{R} \cdot v = \{w + \lambda v \in V_\theta \mid \lambda \in \mathbb{R}\}.$$

For the notion of IRHS we then have the following

**Definition 1** *A* dichotomous item response hypersurface (IRHS) *is a* $D = \dim(V_\theta)$ *dimensional smooth submanifold* $M$ *of* $V_\theta \times [0, 1]$, *so that for any two vectors* $v, w \in V_\theta$ *the intersection of* $(w + \mathbb{R} \cdot v) \times [0, 1]$ *and* $M$ *is a graph of a monotonic function* $w + \mathbb{R} \cdot v \to [0, 1]$.

*We shall say that a* MIRT model *is given when an IRHS is given.*

Note, that while $w + \mathbb{R} \cdot v$ is not canonically isomorphic to $\mathbb{R}$, monotonicity of the map

$$f_{v,w} : w + \mathbb{R} \cdot v \to [0, 1], \quad \lambda \mapsto f_{v,w}(w + \lambda v) \tag{13}$$

can be unambiguously defined by requiring that either $f_{v,w}(w + \lambda v) \leq f_{v,w}(w + \mu v)$ or $f_{v,w}(w + \lambda v) \geq f_{v,w}(w + \mu v)$ for all $\lambda, \mu \in \mathbb{R}$ whenever $\lambda \leq \mu$. We shall use the notation $f_v = f_{v,0}$. Figure 2 shows the intersection of $(w + \mathbb{R} \cdot v) \times [0, 1]$ and $M$ in two dimensions.

To understand the definition better, let us first assume that we choose $v$ to be arbitrary and $w = 0$ in the definition above. Then, the line $w + \mathbb{R} \cdot v = \mathbb{R} \cdot v$ can be understood as an ability direction. The monotonicity requirement of Definition 1 asks for the natural feature that as the ability given by $v$ increases the probability of the correct response increases as well. For non-zero $w$ the requirement is equivalent to the conditional probability of correct response being monotonic with respect to one ability when the rest of the abilities are fixed to a certain not necessarily zero value. To be precise, we should say that for $w \neq 0$ there exists a basis of $V_\theta$ so that the monotonicity requirement reads as the interpretation above. Furthermore, for any basis of $V_\theta$ Definition 1 will ensure the monotonicity of the conditional probability of correct response for any ability direction given any fixed values for the rest of the ability directions (as defined by the basis).

Note also that the collection of maps $f_{v,w}$ for $v, w \in V_\theta$ defines the IRHS completely. For this reason, we shall use the notation $f^M$, or $f$ if no confusion may arise, for the function describing the IRHS $M \subset V_\theta \times [0, 1]$.

One may be tempted to object to the use of notions like manifold and hypersurfaces. It is very important to note, however, that the conditional probability of correct response has been given by a hypersurface in the usual MIRT literature as well. One major difference in terminology is that it was still called *surface* in any dimension, which is a correct usage only in dimension two. In higher dimensions, the object at hand is a hypersurface, a special case of higher dimensional manifolds.

A basis $\mathfrak{v} = (v_1, \ldots, v_D)$ in $V_\theta$ defines a unique isomorphism

$$i_\mathfrak{v} : V_\theta \to \mathbb{R}^D, \quad \sum_{i=1}^{D} \lambda_i v_i \mapsto \sum_{i=1}^{D} \lambda_i e_i, \quad (\lambda_i \in \mathbb{R}), \tag{14}$$

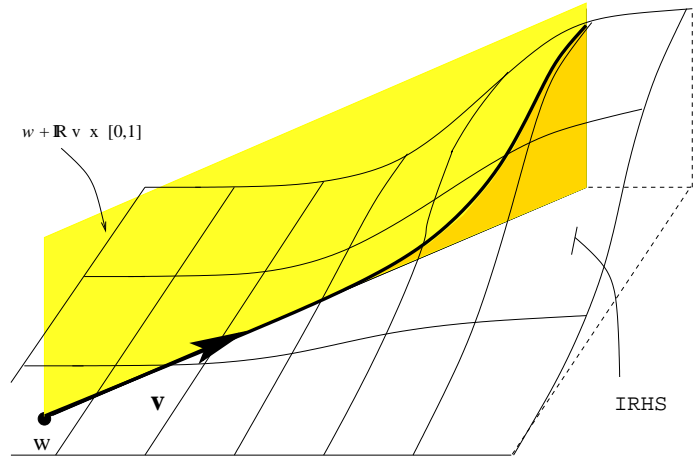

FIG 2. *Intersection of* $w + \mathbb{R} \cdot v \times [0,1]$ *and IRHS. The intersection is the bold curve, which is required to be monotonic.*

where $(e_i)_{i=1}^{D}$ is the standard basis of $\mathbb{R}^D$: $(e_i)_j = \delta_{ij}$ ($\delta_{ij}$ is the Kronecker delta). This isomorphism can be trivially extended to a diffeomorphism

$$\psi_{\mathfrak{v}} : V_{\theta} \times [0,1] \to \mathbb{R}^D \times [0,1], \quad (v,t) \mapsto (i_{\mathfrak{v}}(v), t) \tag{15}$$

and via this diffeomorphism we may transfer the IRHS from $V_{\theta} \times [0,1]$ to $\mathbb{R}^D \times [0,1]$. Now, in $\mathbb{R}^D \times [0,1]$ the image of the IRHS may be given by the graph of a smooth function $f : \mathbb{R}^D \to [0,1]$. Note, however the important difference between using a functional representation like this latter one and using the hypersurface representation directly in $V_{\theta} \times [0,1]$. The functional representation depends on the basis we chose to establish the diffeomorphism $\psi_{\mathfrak{v}}$ and different bases may result in different functional representations.

*Note:* It is tempting to extend this definition to polytomous multidimensional items by defining the *polytomous collection of item response hypersurfaces* for a polytomous item by requiring that the above discussed intersection be a collection of unidimensional polytomous item response curves as produced by some unidimensional polytomous IRT model (e.g. Muraki's partial credit model [9]). The investigation of this possibility is postponed for a forthcoming paper.

### 3.4. Properties of IRHS

In this section we prove the main theorem of the paper. For the sake of transparency, we start with the two dimensional case which is then followed by the more involved general theory.

*3.4.1. Two Dimensional Case*

Using the monotonicity of the model we can prove an interesting elementary property.

**Lemma 1** *In any 2 dimensional MIRT model there exists a line in $V_\theta$ through the origin so that $f_v$ is constant.*

*Proof:* Let us choose a vector $v \in V_\theta$. Note, that if $\inf_{\lambda \in \mathbb{R}} f_v(\lambda v) = \sup_{\lambda \in \mathbb{R}} f_v(\lambda v)$ the lemma is proved, the sought after line is $\mathbb{R} \cdot v$. Therefore, we may assume that $\inf_{\lambda \in \mathbb{R}} f_v(\lambda v) < \sup_{\lambda \in \mathbb{R}} f_v(\lambda v)$. For such a vector either $\lim_{\lambda \to \infty} f_v(\lambda v) = \sup_{\lambda \in \mathbb{R}} f_v(\lambda v)$ or $\lim_{\lambda \to -\infty} f_v(\lambda v) = \sup_{\lambda \in \mathbb{R}} f_v(\lambda v)$. Let $P$ be the set of vectors satisfying the first and $N$ be the set of vectors satisfying the second condition. Both of these sets are non-empty and by continuity, both sets are open. Also, they are clearly disjoint. Therefore, there is a vector $u \in V_\theta$ so that $u \notin N \cup P$. Along the line $\mathbb{R} \cdot u$ the function $f$ is constant. ∎

Note that the proof only uses monotonicity with $w = 0$. Utilizing it for general $w$ the same argument provides the following

**Lemma 2** *In any 2 dimensional MIRT model, through any point $w \in V_\theta$ there exists $v \in V_\theta$ so that along the v-directed line going through $w$ the function $f_{v,w}$ is constant.*

We introduce the term *w-constant line*, or simply *constant line*, for the *v*-directed line going through $w$ as in Lemma 2.

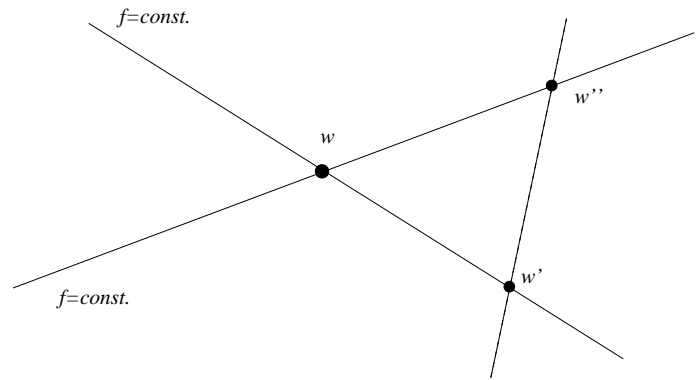

FIG 3. *Non-unique constant line results in constant MIRT model.*

Analyzing the properties of these constant lines further we see that they are actually parallel to one another. That is we have the following

**Lemma 3** *Let $w, w' \in V_\theta$ be two points. Let $v, v' \in V_\theta$ be the corresponding directions of the two constant lines. Then $v = \mu v'$ for some $\mu \in \mathbb{R}$.*

*Proof:* First, we note that if there is a point $w \in V_\theta$ so that there exist two $w$-constant lines, then the model is trivial ($f$ is constant) and the statement is true. For, let $w'$ and $w''$ be the intersections of a general position line in $V_\theta$ with the two $w$-constant lines, respectively. Because $f$ is monotonic along this line, and $f(w') = f(w'')$ $f$ is constant between $w'$ and $w''$, that is $f(tw' + (1-t)w'') = f(w')$ for all $t \in [0, 1]$. Using this argument for every line in general position proves that $f$ is constant everywhere (see Figure 3).

Now, we assume that the constant lines through $w$ and $w'$ are unique. If the two lines are not parallel then they will have an intersection and an argument similar to the previous one shows that $f$ is constant. ∎

The corollary of the previous observation is the

**Theorem 1** *Any 2 dimensional MIRT model is a trivial extension of a unidimensional IRT model.*

*Proof:* We saw in Lemma 3 that a 2 dimensional IRHS is nothing but a collection of parallel lines. Let $v \in V_\theta$ be the direction of these lines. Choosing a transversal $\mathbb{R} \cdot u$ (a line that intersects all of them) to this collection the IRHS can be given by the function $f_u : \mathbb{R} \cdot u \to [0, 1]$. For, let us express an arbitrary $w \in V_\theta$ as a unique linear combination $w = \mu u + \lambda v$ and write

$$f(w) = f(\mu u + \lambda v) = f_u(\mu u). \tag{16}$$

This function $f_u$ can be thought of as a unidimensional IRT model. ∎

### 3.4.2. D-Dimensional Case

Technically, the $D$ dimensional case is not that much more complicated than the 2 dimensional one. It is just much more difficult to visualize the corresponding geometric objects. As we pointed out earlier, the conditional probability "surface" is not 2 dimensional, so strictly speaking it is not a surface in higher dimensions. Our three dimensional training does not allow us to "see" objects in higher dimensions. The formalism we built in the previous section, however, will be applicable, with appropriate modifications, to this situation as well.

The proof of Lemma 1 works for any dimensions. Applying the monotonicity argument for arbitrary $(v, w)$ as above proves the corresponding

**Lemma 4** *In any MIRT model there exists a hyperplane $H_w$ in $V_\theta$ through $w \in V_\theta$ so that $f_{H_w}$ is constant.*

*Proof:* Here, $f_{H_w}$ is the restriction of $f$ to the hyperplane $H_w$. As before, we prove $w = 0$ explicitly; the general case follows the same argument. Let us, as before, define the open sets $P$ and $N$ and note that $P = -N$. Exclude the trivial case of $P = \emptyset$. It is clear that $P \cup N \neq V_\theta$. Locally the boundary of $P$ (the closure of $P$ minus $P$) is a $D - 1$ dimensional submanifold ($D = \dim V_\theta$). Therefore there exists a collection $(c_1, \ldots, c_{D-1})$ of points in $V_\theta \backslash (P \cup N)$ so that $(c_1 - w, \ldots, c_{D-1} - w)$ spans a hyperplane $H_w$. Now, along the line segment joining any point on $\mathbb{R} \cdot (c_i - w)$ with any point on $\mathbb{R} \cdot (c_j - w)$, for some $i \neq j$, the restriction of $f$ should be constant (Figure 3).

Repeating this argument for each pair of line segments shows that along the entire hyperplane $f$ is constant. If there is another hyperplane with this property, then $P = \emptyset$, which is excluded. ∎

Now, a $D - 1$ dimensional hyperplane is to a $D$ dimensional space as a line is to the plane. Using this intuition it is not difficult to adapt the formal proof of Lemma 3 to prove

**Lemma 5** *Let $w, w' \in V_\theta$ be two points. Let $H_w, H_{w'} \subset V_\theta$ be the corresponding two constant hyperplanes. Then $H_w$ and $H_{w'}$ are parallel.*

Now, we are ready to rephrase our main theorem in arbitrary dimension.

**Theorem 2** *Any MIRT model is a trivial extension of a unidimensional IRT model.*

For the sake of explicitness let us write $f^M$ for an arbitrary IRHS in terms of univariate IRT model. Let us fix a transversal $u \in V_\theta$ to the collection of constant hyperplanes. First, we observe that for any $w \in V_\theta$ there is a unique decomposition $w = \mu u + \lambda v$ with $v \in H_w$. Then,

$$f^M(\mu u + \lambda v) = f_u^M(\mu u). \tag{17}$$

Note that if we choose the usual 2Pl or 3PL models the construction yields the Scalar Product model. It is also interesting to note that the MIRT generalization of the Rasch model is equivalent to the generalization of the 2PL model. This is because, while within the univariate Rasch model one may assume that the slope is fixed, when more dimensions are considered simultaneously the assumption of equal slopes is not valid. The relative positions of slopes to one another should be determined during the estimation procedure in lack of *a priori* information.

This kind of models were called generalized compensatory models (GMIRT) in [13]. The link function of an IRHS as GMIRT is $f_u^M$.

### 3.4.3. Absolute Functional Representation for the Scalar Product Model

A notable feature of the Scalar Product model is that using the dual of a vector space it can be defined without referring to coordinates even in its functional form. First, we recall that the *dual $V^*$* of a finite dimensional vector space $V$ is the finite dimensional vector space of the same dimension of linear maps $V \rightarrow \mathbb{R}$:

$$V^* := \{p : V \rightarrow \mathbb{R} \mid p \text{ is linear}\}. \tag{18}$$

The *duality* is the obvious map

$$( \mid ) : V^* \times V \rightarrow \mathbb{R}, \ (p, v) \mapsto ( p \mid v) := p(v). \tag{19}$$

That is, for any $p \in V^*$ and $v \in V$ the quantity $(p \mid v)$ is a real number. It is important to note that the duality, unlike a scalar product, does not involve any choice.

Now, if in MIRT we make the choice, that the ability is modeled by the vector space $V_\theta$ as before and the item is modeled by the discrimination $a \in V_\theta^*$ in the dual space and a real number $b$ then the IRHS of the model is given as the graph of the following function:

$$f_{a,b}^d : V_\theta \to [0,1], \; f_{a,b}^d(\theta) := \frac{1}{1 + e^{-(a \mid \theta) - b}}. \tag{20}$$

In addition to its very satisfying and elegant nature this model has the computational advantage of having the same functional representation in *any* coordinate system. As we shall see later the dimension-wise independent model does not share this nice invariance property.

### 3.4.4. Interpretation of Main Theorem

The statement of the main theorem excludes many existing MIRT models from the pool of monotonic MIRT models. The author's reading of the main theorem is that the only relevant MIRT model is the one defined in (17). This interpretation is backed by the fact the widely used and tested estimation tools exist only for the Scalar Product model, the most relevant of the above extensions ([10]). In the view of Theorem 2 there seems to be a good reason behind that. It seems that lack of monotonicity prevents one to maximize the likelihood function of MIRT models excluded by our approach. This certainly defines a valid future research direction. Also, the existence of an elegant coordinate free functional representation makes the Scalar Product model even more appealing.

On the other hand, model building always has many steps that cannot be entirely backed by theoretical considerations. The process sometimes is dictated by personal preferences and tastes. It is possible that some readers may not be willing to except the requirement of monotonicity as formulated in Definition 1 as a crucial and necessary feature of an MIRT model. For those readers the main theorem is interpreted a bit differently. First, we note the close connection between the notion of *compensatory* model to monotonicity. Usual terminology is that the model is compensatory, if the probability of the correct response may be high even with the lack of ability in all but one dimension. That is, sufficiently high ability in one dimension is able to compensate for the lack of it in other dimensions. In fact, compensatory property follows from monotonicity as an easy application of Theorem 2. If compensatory property is understood in a sense that it is true in any coordinate system, then the reverse is also true, and the two notions are equivalent. With this in mind the theorem states that *any compensatory MIRT model is a direct generalization of a univariate IRT model.*

In either way, Theorem 2 establishes a prominent role for the Scalar Product model as an MIRT model.

### 3.5. Estimation in MIRT

Let us now restrict our attention to the Scalar Product model. A typical two dimensional ($D = 2$) student likelihood (12) is given in Figure 4. As in unidimensional IRT, the maximum place of this function plays a special role in the estimation of MIRT model parameters. A curious feature of this graph is that a pronounced unbalance can be observed between the standard errors of the two ability estimates. Here, standard error is understood as the inverse of the curvature of the graph at the maximum place. There is a well identified direction in which the standard error is minimal and in the direction orthogonal to this the standard error appears to be much bigger. One may even say that, despite our efforts, the model shows definite signs of unidimensionality.

The reason behind this is very simple. A student likelihood is formed as a product of probabilities of the actual responses given by item response hypersurfaces similar to the one shown on the RHS of Figure 1. These hypersurfaces are always increasing towards the first quadrant (correct response) or towards the third quadrant (incorrect response). Hence, the product of these will be the above observed "ridge" of Figure 4. It is a ridge because the observed response is either correct or incorrect and no distinction is made between events of the students using only one of the dimensions correctly during the assessment. In other words, since there is no observed data for the different dimensions, the model will not be able to provide two distinct, meaningful estimates for the abilities of the person on the different cognitive dimensions.

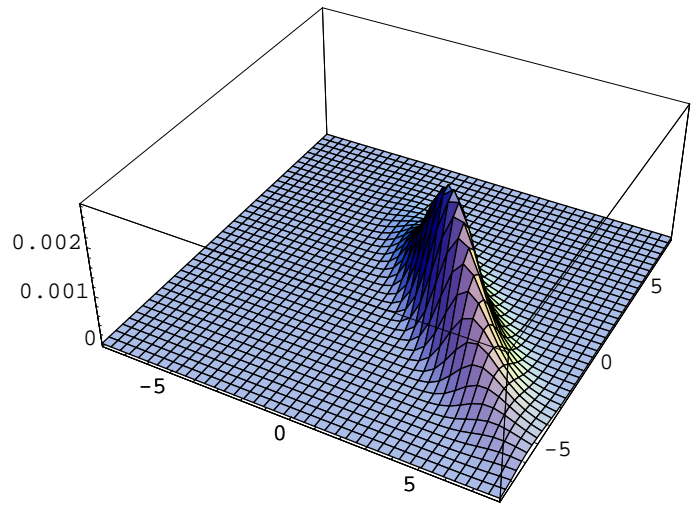

Fig 4. *Scalar Product MIRT student likelihood.*

### 3.6. Dimension-wise Independence

The careful reader should have noticed, that concerning one particular point the presentation is not faithful to its own principles. That is, the notion of dimension-wise independence was used without any discussion of its invariance, or coordinate system independence. It is easy to see that the dimension-wise independent model does not satisfy the requirement of monotonicity, therefore we would not consider it as a valid MIRT model. On the other hand, it might be useful to see explicitly how badly the the functional representation of the dimension-wise independent model behaves to appreciate the niceties of the Scalar Product model even more.

Invariance of dimension-wise independence for the the model

$$f_{a,b}^w : \mathbb{R}^D \to [0,1], \quad f_{a,b}^w(\theta) = \prod_{d=1}^{D} \frac{1}{1 + e^{-a_d(\theta_d - b_d)}} \tag{21}$$

would mean that the factorization property holds in any other coordinate systems.

Mathematically, this would require that for any invertible matrix $G \in \mathrm{GL}(D)$ ($G$ expresses change of coordinates) we have a function

$$h^G : \mathbb{R} \times \mathbb{R} \times \mathbb{R} \to \mathbb{R}, \quad (a_1, b_1, t) \mapsto h^G(a_1, b_1, t)$$

and a pair of invertible matrices $U, V \in \mathrm{GL}(D)$ so that when $\theta = G \cdot \theta'$ ($\theta, \theta' \in \mathbb{R}^D$) we have a factorization

$$f_{a,b}^w(\theta) = \prod_{d=1}^{D} h^G((Ua)_d, (Vb)_d, \theta'_d) \tag{22}$$

that is

$$
\begin{aligned}
f_{a,b}^w(\theta) &= f_{a,b}^w(G \cdot \theta') \\
&= \prod_{d=1}^{D} \frac{1}{1 + e^{-a_d(\sum_{d'=1}^{D} g_{dd'} \theta'_{d'} - b_d)}} \\
&= \prod_{d=1}^{D} h^G(a'_d, b'_d, \theta'_d),
\end{aligned}
\tag{23}
$$

with $a'_d = (Ua)_d$ and $b'_d = (Vb)_d$. The role of $U$ and $V$ is to ensure that the function $h^G$ is the same for all factors in the product by allowing this function to depend on different linear combinations of the elements of $a$ and of $b$.

To show that this is too much to ask for in general, let us first assume that a factorization $f(x, y) = h(x)g(y)$ holds for some function $f$ so that $h(0) \neq 0$ and

$g(0) \neq 0$. Then,

$$h(x) = \frac{f(x,0)}{g(0)},$$

$$g(y) = \frac{f(0,y)}{h(0)},$$

and

$$\frac{1}{h(0)g(0)} = \frac{f(x,y)}{f(x,0)f(0,y)}. \tag{24}$$

Now, for the sake of concreteness, let us take $D = 2$ and $a = (a_1, a_1) \in \mathbb{R}^2$ and $b = (0,0) \in \mathbb{R}^2$. Also, let us take $G = \begin{pmatrix} 1 & 1 \\ 1 & -1 \end{pmatrix}$. With these, (23) becomes

$$\frac{1}{1 + e^{-a_1(\theta_1' + \theta_2')}} \cdot \frac{1}{1 + e^{-a_1(\theta_1' - \theta_2')}} = h(\theta_1')g(\theta_2') \tag{25}$$

with some $h, g : \mathbb{R} \to \mathbb{R}$. From (24) the function

$$\frac{1}{h(0)g(0)} = \frac{\frac{1}{1 + e^{-a_1(\theta_1' + \theta_2')}} \cdot \frac{1}{1 + e^{-a_1(\theta_1' - \theta_2')}}}{\frac{1}{(1 + e^{-a_1\theta_1'})^2} \cdot \frac{1}{1 + e^{-a_1\theta_2'}} \cdot \frac{1}{1 + e^{+a_1\theta_2'}}} \tag{26}$$

should be constant. This is clearly not the case, showing that the factorization (23) does not hold in general.

It seems that the definition of dimension-wise independence is not an absolute one. We have a choice of either dropping it altogether, or if need arises, we may change it. To formulate this notion we have to relax the monotonicity requirement of MIRT in Definition 1 by requiring the monotonicity of the $f_v$ for all $v \in V_\theta$, that is assumed $w$ is zero in Definition 1 and in 13. Let us call this type of models *ray-wise monotonic MIRT* models.

**Definition 2** *A ray-wise monotonic MIRT model given by an IRHS is* dimension-wise independent *if there exists a coordinatization of abilities so that the functional representation of the model $f_{a,b}(\theta)$ can be written as a product of factors*

$$f_{a,b}(\theta) = \prod_{d=1}^{D} h(a_d, b_d, \theta_d'). \tag{27}$$

The specialty of this property comes from the fact that for a general IRHS it is very rare that the functional representation can be factored so that one may consider it dimension-wise independent. This interpretation was used throughout the paper, when the Whitley model was called dimension-wise independent.

## 4. Conclusion

A coordinate free definition of MIRT has been put forward in the paper. Our main argument is that in a coordinate free setup it is easier to tell apart genuine

MIRT objects from potential artifacts. These artifacts can be notions and relationships that should not be considered integral parts of the model since their key features which may be apparent in one could vanish in another coordinate system. We showed that it is possible to provide a full classification of monotonic models solely based on general, coordinate-free considerations.

It is important to note that the classification was carried out at the level of a single item IRHS, but it is in no way restricted to single item tests. IRT and MIRT models handle tests by invoking the local independence assumption and form the likelihood of the model by multiplying single item conditional probabilities together. The *flavor* of the test (Rasch, 2Pl, normal ogive, compensatory, polytomous, etc.) is always given at the single item level. Our treatment is no exception.

It is very important that the reader does not mistake the promotion of the coordinate free description for an argument for a completely coordinate free handling of the entirety of MIRT. In fact, it should be explicitly stated that without a choice of coordinates meaningful MIRT practice cannot exist. In addition to this, every discussion of MIRT features can be fully carried out using $\mathbb{R}^D$ as the main model space for abilities. Should such a path be chosen, however, one has to be careful to meticulously maintain the coordinate system invariance of the theory every step of the way. The contribution of this paper is an introduction of a framework to ease this burden by keeping the presentation *absolute* (without choosing any coordinates) for as long as possible. The paper shows that one may be able to formulate general statements and reach valuable insights before switching to *relative* mode by an introduction of a particular basis. It is likely that someone may observe the relevance of a notion while in a particular coordinate system and may want to establish whether it is invariant by trying to create a definition in the absolute framework presented here.

It is noteworthy, that the necessity of the existence of a coordinate free representation of our physical world led Einstein to formulate both the special and the general theories of relativity ([2; 3]). The fundamental dogma in relativity theory is that the events of the physical world take place without being aware of any coordinate system. Therefore, any faithful description should be invariant of the change of coordinate system. Better yet, a description of the physical world is sought that bypasses the use of coordinates altogether.

A reader interested in the successes of coordinate free description of the physical world may also find the books [6; 7] useful.

## Acknowledgments

The author would like to extend his gratitude to Paul Holland for his the insights and constant encouragement. The author is also indebted to Shelby Haberman and Henry Chen for fruitful discussions. This paper was inspired by the thought provoking presentation of Mark Reckase at the *Educational Testing Service* in September, 2006. The author is indebted for the lively discussion during this presentation.

## References

[1] A. Birnbaum. Some latent trait models and their use in inferring an examinee's ability. In F. M. Lord and M. R. Novick, editors, *Statistical Theories of Mental Test Scores*, pages 397–479. Reading, MA: MIT Press, 1968.

[2] Albert Einstein. Zur Elektrodynamik bewegter Körper. *Ann. Phys.*, 17:891–921, 1905.

[3] Albert Einstein. Grundlagen der allgemeinene Relativitätstheorie (The foundation of the general theory of relativity). *Ann. Phys.*, 49(4):284–339, 1916.

[4] Gerhard H. Fischer. On the existence and uniqueness of maximumlikelihood estimates in the Rasch model. *Psychometrika*, 46(1):59–77, 1980. MR0655008

[5] Paul R. Halmos. *Finite-Dimensional Vector Spaces*. New York: Springer, 2 edition, 1974. MR0409503

[6] Tamás Matolcsi. *A Concept of Mathematical Physics: Models in Mechanics*. Akadémiai Kiadó, Budapest, 1986. MR0873263

[7] Tamás Matolcsi. *Spacetime without Reference Frames*. Akadémiai Kiadó, Budapest, 1993. MR1240055

[8] R. L. McKinley and M. D. Reckase. The use of the general Rasch model with multidimensional item response data. *Research Report ONR 82-1*, 1982.

[9] Eiji Muraki. A generalized partial credit model: Application of an EM algorithm. *Appl. Psychol. Meas.*, 16:159–176, 1992.

[10] M. D. Reckase. A linear logistic multidimensional model for dichotomous item response data. In W. J. van der Linden and R. K. Hambleton, editors, *Handbook of Modern Item Response Theory*, pages 271–286. New York: Springer, 1997.

[11] F. W. Warner. *Foundations of Differentiable Manifolds and Lie Groups*. Scott, Foresman and Company, Glenview, Illinois, 1971. MR0295244

[12] Susean E. Whitely. Measuring aptitude processes with multicomponent latent trait models. *Technical Report No. NIE-80-5*, 1980.

[13] Jinming Zhang and William F. Stout. Conditional covariance structure of generalized compensatory multidimensional items. *Psychometrika*, 64:129–152, 1999. MR1700706